\documentclass[11pt]{article}
\usepackage{amsmath, amssymb, amsthm}
\usepackage{geometry}
\usepackage{url, hyperref}
\usepackage[utf8]{inputenc}
\usepackage[french, english]{babel}
\usepackage{comment}

\newtheorem{THEOREM}{Theorem}[section]

\newtheorem{theorem}[THEOREM]{Theorem}
\newtheorem{lemma}[THEOREM]{Lemma}
\theoremstyle{definition}
\newtheorem{definition}[THEOREM]{Definition}
\newtheorem{claim}[THEOREM]{Claim}

\newtheorem{example}[THEOREM]{Example}
\newtheorem{fact}[THEOREM]{Fact}

\newtheorem{Observation}[THEOREM]{Observation}
\newenvironment{observation}{\begin{Observation}}
{\end{Observation}}
\newtheorem{Notation}[THEOREM]{Notation}

\newcommand{\forces}{\Vdash}
\newcommand{\dom}{\text{dom}}

\newcommand{\ran}{\text{ran}}

\newcommand{\M}{\mathbf M}

\newcount\skewfactor
\def\mathunderaccent#1#2 {\let\theaccent#1\skewfactor#2
\mathpalette\putaccentunder}
\def\putaccentunder#1#2{\oalign{$#1#2$\crcr\hidewidth
\vbox to.2ex{\hbox{$#1\skew\skewfactor\theaccent{}$}\vss}\hidewidth}}
\def\name{\mathunderaccent\tilde-3 }

\newcommand{\rest}{\upharpoonright}

\newcommand{\cf}{\text{cf}}

\newcommand{\eop}{\bigstar}
\newcommand{\DD}{\mathcal D}
\newcommand{\EE}{\mathcal E}
\newcommand{\FF}{\mathcal F}

\newcommand{\eod}{\scriptstyle\blacksquare}

\title{Adding unbounded and not bounded sets to a singular cardinal}

\author{Mirna D{\v z}amonja\\
Logique Consult, Liège, Belgium}

\date{}

\begin{document}

\maketitle\footnote{Earlier versions of this paper were co-authored with Vladimir Kanovei, who had to stop the collaboration due to a potential conflict of interest.}

\begin{abstract} This paper solves a problem (SP) attributed to Solovay in the 1960s, which is to
add an unbounded sequence to a singular cardinal $\kappa$ without adding bounded subsets to $\kappa$, while preserving 
$\kappa^+$. The solution applies to any singular cardinal and any cofinality of the singular cardinal, with the assumption that $\kappa$ is a strong limit needed for the case of uncountable cofinality. In particular, the solution applies to $\kappa=\aleph_\omega$.
Unsuccessful attempts to solve (SP), mostly for the countable cofinality, have appeared in published papers by various authors. 
\end{abstract}

\section{Introduction and Motivation}  In their 2020 paper \cite{DobrinenHathawayPrikry2020} Dobrinen, Hathaway and Prikry present a method which they attribute to an unpublished paper by Prikry in the late 1960s. Its main purpose is to solve a problem posed by Solovay in the 1960s, which is:

{\bf Solovay's problem (SP)}: Find a forcing to add an unbounded sequence to a singular cardinal $\kappa$ without adding any bounded sequences to $\kappa$, and while preserving $\kappa^+$.

Unfortunately, we are not convinced by the arguments from \cite{DobrinenHathawayPrikry2020}, as the literal translation of the ``usual arguments" from arboreal forcings on $\omega$ which they relay upon, fails at many points in their intended argument. One of them is the compactness argument needed for certain topological arguments about $2^\omega$ and not present in their space $N$. A simplification of the argument from  \cite{DobrinenHathawayPrikry2020} was intended in \cite{zbMATH07903059}, but just like \cite{DobrinenHathawayPrikry2020}, it fails an attempted fusion argument\footnote{In an e-mail dated 04 September 2026, Heike Mildenberger claimed that she had an explanation to the argument from \cite{zbMATH07903059} and she reiterated it in an e-mail dated 10 November 2026. Unfortunately, it enforces the incorrect attempt of using a Cantor-Bendixon argument to address the question of not adding bounded subsets to $\kappa$, in a situation where this argument does not apply.}. These papers have a predecessor \cite[\S 1.3]{GitikPrikrylike} which claims to have obtained the result from an $\omega$-sequence of measurables. However, it also fails an attempted fusion argument. The difficult part in all these approaches is to obtain that no bounded sequences are added to the singular cardinal. The purpose of this paper is to define a family of forcing notions ${\mathbb P}_{\bar{\kappa}}[\bar{\FF}]$ and prove

\begin{theorem}[{\bf The Solution of (SP)}]\label{the:main-result-unc} Suppose that $\kappa$ is any singular cardinal of cofinality 
$\aleph_0\le \theta<\kappa$ and that $\bar{\kappa}$ is an increasing $\theta$-sequence of regular cardinals
with supremum $\kappa$. In addition, if $\theta>\aleph_0$, assume that $\kappa$ is a strong limit. 

Then the forcing 
${\mathbb P}_{\bar{\kappa}}[\bar{\FF}]$ adds a $\theta$-sequence to $\kappa$ without adding any bounded subsets to $\kappa$, and it preserves 
$\kappa^+$.

(Note that the theorem applies to the situation $\kappa=\aleph_\omega$.)
\end{theorem} 

The problem of adding a subset to a singular cardinal without adding any bounded subsets to it was addressed by a different method in 1973 by Kanovei, \cite{Kanovei1973}. Forcing over a model of $\mathbf{V}=\mathbf{L}$, for any cardinal $\lambda$, he constructs a forcing extension in which no new subsets of size $<\lambda$ are added to any ordinal $(<\lambda^+)^{\mathbf{L}}$, while a subset is added to $\lambda$ (\cite[Th 2]{Kanovei1973}). This does not solve (SP), because it collapses $\lambda^+$. 

We address (SP) by a combination of ideas from \cite{FriedmanGitman:ModelOfACNotDCInaccessible}, \cite{Property-B} and some new ideas.
To explain this, let us first concentrate on the situation of $\cf(\kappa)=\aleph_0$. In the first place we adapt the Generalised Perfect Set Forcing from \cite{Property-B} to a singular cardinal $\kappa$ of countable cofinality, using a sequence 
$\bar{\kappa}=\langle \kappa_n:\,n<\omega\rangle$ of regular cardinals strictly increasing to $\kappa$ and a sequence
$\bar{\FF}=\langle \FF_n:\,n<\omega\rangle$ of filters where each $\FF_{n+1}$ is a $(<\kappa_n)$-complete filter over $\kappa_{n+1}$. This could be 
the trivial filter $\{\kappa_{n+1}\}$- so no large cardinal is needed- or it could be a measure on $\kappa_{n+1}$, if one wants to work with measurables. 

The forcing, which we denote ${\mathbb P}_{\bar{\kappa}}[\bar{\FF}]$, uses trees of height $\omega$ where the splittings on the level $n$ are
into a set in $\FF_n$ and along any branch the set of splitting nodes is infinite. This allows to show that ${\mathbb P}_{\bar{\kappa}}[\bar{\FF}]$ adds an 
$\omega$-sequence to $\kappa$ and preserves $\kappa^+$. The more difficult part of the argument shows that the assumption of the filters being sufficiently closed gives that the forcing also does not add bounded sequences to $\kappa$. Hence we obtain a solution to the problem (SP) for the countable cofinality.

The situation of adding subsets to singulars of countable cofinality, although not the problem (SP) directly, was considered by the papers \cite{zbMATH06923337} and \cite{zbMATH07432945}. The paper \cite{zbMATH06923337} introduces a forcing called AIM forcing, which deals with the situation of a singular $\kappa$ of countable cofinality which is a limit of $\lambda$-supercompact cardinals for some $\lambda$ inaccessible $>\kappa$. It is shown in the main theorem of
the paper that AIM forcing does not add bounded subsets of $\kappa$. However the forcing does not solve (SP) even under these large cardinal assumptions, since the theorem also shows that $\kappa^+$ is collapsed, and in fact $\kappa^+$ of the extension is $\lambda$.
In \cite{zbMATH06923337} it is shown that this forcing fits the general framework of $\Sigma$-Prikry forcing introduced in that paper. It is possible that some other example of $\Sigma$-Prikry forcing might be applicable to (SP) of countable cofinality under large cardinals, although none of the presently known 
results fit this description. The method of $\Sigma$-Prikry forcing has some similarities to ours, in that it considers layering forcings into countably many suborders and that it works with a sequence of filters, but our forcing does not fit the axiom scheme of $\Sigma$-Prikry forcing directly.

Let us now describe the situation of arbitrary cofinality of the singular. Levine and Mildenberger \cite[Th. 3.1]{zbMATH07903059}, who have claimed to have simplified the proof from  \cite{DobrinenHathawayPrikry2020}, also proved that the two forcings from their paper, generalisations of Laver and Miller forcing to this context, do not solve (SP) for arbitrary cofinality. Our method does.

For this result we need to modify the definition
of ${\mathbb P}_{\bar{\kappa}}[\bar{\FF}]$, which we still call by the same name, since the version of ${\mathbb P}_{\bar{\kappa}}[\bar{\FF}]$ used for the 
proof of countable cofinality is a special case of the more general version. This modification is done by blending the definition of 
${\mathbb P}_{\bar{\kappa}}[\bar{\FF}]$ from the proof of the countable case with the Generalised Perfect Set Forcing  from \cite{Property-B}. 

In particular, Theorem \ref{the:main-result-unc} shows that (SP) does not require large cardinals\footnote{In an e-mail dated 07 September 2026, Moti Gitik claimed that the argument for the countable cofinality would need an $\omega$-sequence of measurable cardinals. In his e-mail dated 11 November 2026, Mohammad Golshani claimed that a measurable cardinal is needed. In our opinion, both claims have a problem when using the covering lemma.}. \smallskip

The paper is organised as follows. In \S\ref{sec:the-original-Prikry-modification}  we explain what we think that 
\cite{DobrinenHathawayPrikry2020} might have meant for the case of countable cofinality and prove that our version adds an unbounded $\omega$-sequence to $\kappa$ and preserves $\kappa^+$.
In \S\ref{sec:the-original-Prikry-modification} we introduce a family of related forcing notions
${\mathbb P}_{\bar{\kappa}}[\bar{\FF}]$, of which the forcing from \S1 is a special case. We show that none of them adds a bounded subset to $\kappa$.

In \S\ref{sec:with-filters} we recall the Generalised Perfect Set Forcing from 
\cite{Property-B} and in \S\ref{sec:unc-cof} we use it to define a version of the forcing ${\mathbb P}_{\bar{\kappa}}[\bar{\FF}]$ which works for any cofinality of the singular cardinal. Hence we prove Theorem \ref{the:main-result-unc} and solve (SP).

{\bf Acknowledgement}. I would like to gratefully acknowledge Mohammad Golshani who pointed out a mistake in a previous version of this paper, in the proof of Lemma \ref{lem:no-bounded}. His suggested solution was to assume that cardinals $\kappa_n$ for $1\le n<\omega$ were measurable, but this is not necessary. Clarifying the situation regarding his comment allowed for a more general theorem with a simpler proof. It also allowed us to revisit and correct the forcing possibly intended by Prikry in \cite{DobrinenHathawayPrikry2020}. 

\section{A Forcing Attributed to Prikry}\label{sec:the-original-Prikry-modification} We present the forcing which might have been intended in Dobrinen et al. \cite{DobrinenHathawayPrikry2020}, who attribute that part of their paper to Prikry.\footnote{The authors of \cite{DobrinenHathawayPrikry2020} work with a strong limit cardinal even in the case $\cf(\kappa)=\aleph_0$, but we have not found the use for it.} The approach taken here is to make the forcing as close as possible to the ordinary Perfect Set Forcing on 
$\omega$. 

Suppose that $\kappa$ is a singular cardinal with $\cf(\kappa)=\aleph_0$. Let
$\bar{\kappa}=\langle \kappa_n:\,n<\omega\rangle$ be an increasing sequence of regular cardinals with $\kappa_0=0$, $\kappa_1\ge \aleph_0$ and
$\sup_{n<\omega}\kappa_n=\kappa$. 
We shall work with subtrees of the structure
\[
T^\ast=\bigcup_{n<\omega} (\{n\}\times {}^{<\omega}\kappa_n),
\]
where the elements of $T^\ast$ are partially ordered by the following relation:
\[
(n, s)<^\ast (m,t)
\]
if $n<m$ and $s\subseteq t$ or $n=0$ and $s=\langle \rangle$ is the empty sequence.
Hence, $T^\ast$ is a tree of height $\omega$ and cardinaliy $\kappa$. Notice that its definition depends on $\bar{\kappa}$, so we should in fact denote it by $T^\ast _{\bar{\kappa}}$. We omit ${}_{\bar{\kappa}}$ if the ${\bar{\kappa}}$ we mean is clear from the context.

\begin{definition}\label{def:perfect} Suppose that $T\subseteq T^\ast$ is a tree.

(1) We say that a node $(n,s)\in T$ is {\em terminal} if $|\{t:\;(n,s)<_T(n+1, t)\}|= 0$, otherwise the node is {\em non-terminal}. The non-terminal node is  {\em splitting in} $T$ if $|\{t:\;(n,s)<_T(n+1, t)\}|> 1$, otherwise the non-terminal node $(n,s)$ is a non-splitting node.
A splitting node $(n,s)$ is {\em fully-splitting} if $\{t:\;(n,s)<_T(n+1, t)\}=\kappa_{n+1}$.

{\noindent (2)} We say that a tree $T\subseteq T^\ast$ is {\em $\bar{\kappa}$-perfect} if
\begin{itemize}
\item every node in $T$ is either non-splitting or fully-splitting,
\item every node in $T$ has a fully-splitting extension.
\end{itemize}
(In particular, $T$ has no terminal nodes).

{\noindent (3)} The forcing notion $\mathbb P=\mathbb P_{\bar{\kappa}}$ consists of all $\bar{\kappa}$-perfect subtrees of $T^\ast _{\bar{\kappa}}$ ordered by reverse inclusion.

{\noindent (4)} Suppose that $p$ is a $\bar{\kappa}$-perfect tree and $(o,t)\in p$. We denote by $\le_p$ the order in $p$ and define
\[
\deg_p((o,t))=|\{(m,s):\,(m,s)<_p (o,t) \mbox{ is splitting}\}|.
\]

{\noindent (5)} For $p,p' \in \mathbb P$, we let $p\le_n p'$ if $p\le p'$ and for every $(o,t)\in p$,
\[
\deg_p((o,t))\le n\implies (o,t)\in p'.
\]
$\eod_{\mbox{\tiny{Def.}} \ref{def:perfect}}$
\end{definition}

\begin{observation}\label{obs:decreasing} The relations $\le_n$ are partial orders such that for $m\le n$ we have  $\le_{n}\subseteq \le_m$.
\end{observation}

\begin{definition}\label{def:stem} (1) For a condition $p\in \mathbb P$ we define {\em the stem} $s[p]$ as the unique splitting node of $p$ which has no splitting 
predecessors in $p$. 

{\noindent (2)} For $(n,s)\in p$ we denote by $p_{(n,s)}$  the subtree of $p$ given by
\[
p_{(n,s)}=\{(m,t)\in p:\,(m,t)\le_p (n,s) \vee (n,s) \le_p (m,t)\}.
\]

{\noindent (3)}
If $(n,s)\in p$ we say that $n$ is the {\em length of} $(n,s)$ and we denote $n=\lg((n,s))$. For $k\le n$ we define $(n,s)\rest k$ as the unique $(k,t)\in p$
with $(k,t)\le_p (n,s)$. 

{\noindent (4)} If $s[p]=(n,t)$ then by $s'[p](k)$ for $k<n$ we denote the value of $t(k)$.
$\eod_{\mbox{\tiny{Def.}} \ref{def:stem}}$
\end{definition}

\begin{lemma}[{\bf The generic function}]\label{lem:generic} Suppose that $G$ is a $\mathbb P_{\bar{\kappa}}$-generic filter over the ground model $\mathbf M$. Let us attempt to define a function $f:\,\omega\to\kappa$ by 

$f(n)=\alpha$ iff
there is $p\in G$ with $\lg(s[p])\ge n+1$ and $s'[p](n)=\alpha$.

Then $f:\,\omega\to\kappa$ is a well-defined total function, it is not an element of $M$ and $\sup(\ran(f))=\kappa$.
\end{lemma}

\begin{proof} Let us first see that $f$ is a well-defined function on its domain. Suppose that $p,q\in G$ satisfy that $\lg(s[p])\ge n+1$
and $\lg(s[q])\ge n+1$. Since $G$ is a filter,
the conditions $p$ and $q$ are compatible, so $s[p]\le^\ast  s[q]$ or vice versa. Therefore $s'[p](n)=s'[q](n)$.

Hence, if $n\in \dom(f)$, then $f(n)$ is well-defined.

We exhibit certain dense subsets of $\mathbb P$ which are elements of $\M$. For $n<\omega$ and $\zeta<\kappa$ let 
\[
{\mathcal D}_{n, \zeta}=\{q\in \mathbb P:\,\lg([s[q]) >n \wedge (\exists m\in [n, \lg([s[q]) \, [s'(m)>\zeta]\}.
\]
We show that ${\mathcal D}_{n, \zeta}$ is a dense set in $\mathbb P_{\bar{\kappa}}$. Suppose that $p\in \mathbb P_{\bar{\kappa}}$ is given. Let $m$ be the first such that $\zeta<\kappa_m$ (so $m\ge 1$) and let $(o, t)$ be a splitting node in $p$ with $o\ge m-1$. Therefore, for any extension  $(o+1, t')$ of $(o, t)$ in $p$,
we have that $t$ is a function whose values are in $\kappa_{o+1}$ and all ordinals in $\kappa_{o+1}$ are taken as values of
such $t'(o)$. Hence there is $(o+1, t')\in p$ with $t'(o+1)>\zeta$.
Then letting $q=q_{(o+1,t')}$ gives
$q\ge p$ and $q\in {\mathcal D}_{n, \zeta}$.
This argument shows that ${\mathcal D}_{n, \zeta}$ is dense. Any $p\in G\cap {\mathcal D}_{n, \zeta}$ witnesses that $f(n)$ is well defined and that
$\sup(\ran(f))\ge\zeta$. Since $n<\omega$ and $\zeta<\kappa$ were arbitrary, the genericity of $G$ over $\M$ implies that $f$ is a total function whose
range is unbounded in $\kappa$.

For $g:\,\omega\to \kappa$ with $g\in \M$ we let 
\[
\EE_g= \{q \in \mathbb P:\, (\exists m<\lg(s[q])\,[s'[q](m)\neq g(m)]\}.
\]
Then $\EE_g\in \M$ and a proof similar to the above one gives that $\EE_g$ is a dense subset of $\mathbb P$. Namely, given $p\in \mathbb P$, we can find a fully-splitting node $(o,t)\in p$ of height at least $1$ and then choose an extension $(o+1, t')$ of $(o,t)$ in $p$ with $t'(o)\neq g(o)$. Then defining $q=q_{(o+1,t')}$ gives $q\ge p$ and $q\in \EE_g$.
By the definition of the sets $\EE_g$ and the genericity of $G$ over $\M$ we have that $f$ is a function which is not in $\M$.
$\eop_{\ref{lem:generic}}$
\end{proof}

\begin{definition}\label{def:gen-sequence} With the notation of Lemma \ref{lem:generic}, we denote $f=g[G]$.
$\eod_{\mbox{\tiny{Def.}} \ref{def:gen-sequence}}$
\end{definition}

\begin{lemma}[{\bf Fusion Lemma}]\label{lem:fusion} Suppose that $\langle p_n:\,n<\omega\rangle$ is a sequence of conditions in $\mathbb P_{\bar{\kappa}}$ such that
$p_n\le_{n} p_{n+1}$ for every $n<\omega$. Then $q=\bigcap_{n<\omega}p_n$ is a condition in $\mathbb P_{\bar{\kappa}}$ and it satisfies that 
for every $n<\omega$ we have $p_n\le_n q$.
\end{lemma}

\begin{proof} Suppose that $m<\omega$ and $(o,t)\in p_m$ is such that $\deg_{p_m}((o,t))\le m$. By Observation \ref{obs:decreasing} we have that $(o,t)\in p_n$ for every $n\ge m$. In particular $(o,t)\in q$. So, $p_m\le_m q$.

To finish the proof, we need to prove that $q$ is a perfect subtree of $T^\ast$. Let $(o,t)\in q$ and let $n^\ast=\deg_q((o,t))$. We observe that the fact that 
$\langle p_n:\,n<\omega\rangle$ is increasing in $\mathbb P_{\bar{\kappa}}$ implies that $\langle \deg_{p_n}((o,t)):\,n<\omega\rangle$ a decreasing sequence
of natural numbers with limit $n^\ast$, and hence there is $n$ such that for every $m\ge n$ we have that $\deg_{p_m}((o,t))=n^\ast$. Let us choose such $n$
so that $n \ge n^\ast$.

Let $(o',t')$ be a fully-splitting extension of $(o,t)$ in $p_{n+1}$ and of the shortest height. In particular, $\deg_{p_{n+1}}((o'+1, t''))=n^\ast+1$, for any $(o'+1, t'')\ge_{p_{n+1}} (o',t')$. By the definition of $\le_{n+1}$ and the fact that $n^\ast+1\le n+1$ we have that $(o', t'), (o'+1,t'')\in q$ for any such $t''$. Hence, $(o',t')$ is a fully-splitting extension of $(o,t)$ in $q$. Therefore, every node in $q$ has a fully-splitting extension.

A similar argument proves that every splitting node in $q$ is fully-splitting. Namely, suppose that $(o,t)$ is a splitting node in $q$ and let $n^\ast=\deg_q((o,t))$.
Then in particular $(o,t)$ is a splitting node in $p_{n+1}$ for all large enough $n$ and all the $p_{n+1}$-immediate successors have degree $n^\ast+1$ in such $p_{n+1}$. Therefore, they continue being immediate successors of $(o,t)$ in $q$, and since $(o,t)$ is fully-splitting in $p_{n+1}$, it continues being fully splitting in $q$.
$\eop_{\ref{lem:fusion}}$
\end{proof}

The following property \ref{Miller-density}(1) was called {\em Mastering $\kappa$} in \cite{Property-B}. It is typical of the variants of Perfect Set Forcing and is used for the
preservation of $\kappa^+$. 

\begin{lemma}[{\bf Mastering $\kappa$ and Preservation of $\kappa^+$ Lemma}]\label{Miller-density} (1) Suppose that $n<\omega$ and $p\forces``\name{\tau}\mbox{ is an ordinal}"$.
Then there is $r\ge_{n+1} p$ and a set $X$ of size $\le\kappa$ such that $r\forces``\name{\tau}\in X"$.

\smallskip 

{\noindent (2)} Forcing with $\mathbb P_{\bar{\kappa}}$ preserves $\kappa^+$.
\end{lemma}

\begin{proof} 
(1) Let $\{s_\alpha:\,\alpha\le \alpha^\ast\}$ be an enumeration of all $s$ with $(m,s)\in p$ for some $m$ and such that and $\deg_p ((m,s))=n$.
So $|\alpha^\ast|\le |p|\le \kappa$. For each $s_\alpha$ let $m_\alpha<\omega$ be such that $(m_\alpha, s_\alpha)\in p$.
For each $\alpha$ and $i<\kappa_{m_\alpha+1}$ such that $(m_\alpha+1, s_\alpha\frown i)\in p$, let $r_{\alpha, i}\ge p_{(m_\alpha+1, s_\alpha\frown i)}$ be some condition that decides the value of $\name{\tau}$, say
$\gamma_{\alpha, i}$. Let $X=\{\gamma_{\alpha, i}:\,\alpha<\alpha^\ast, i<\kappa\}$, so $|X|\le\kappa$.

Now let $r=\bigcup_{\alpha<\alpha^\ast, i<\kappa} r_{\alpha, i}$ and observe that $r\in \mathbb P_{\bar{\kappa}}$ and $r\ge_{n+1} p$. Suppose that $q\ge r$
and $q$ decides the value of $\name{\tau}$. Then $q$ must be compatible with some $r_{\alpha, i}$ and hence the value decided by $q$ for $\name{\tau}$
must be one among the $\gamma_{\alpha, i}$'s. Therefore $r\forces``\name{\tau}\in X"$.

{\noindent (2)} Suppose for a contradiction that $\kappa^+$ is collapsed to $\kappa$ when forcing with $\mathbb P_{\bar{\kappa}}$. In particular, the forcing provides a cofinal $\omega$-sequence in $(\kappa^+)^{\M}$. Let $\name{h}$ be a $\mathbb P_{\bar{\kappa}}$-name for such a cofinal function from 
$\omega$ to $(\kappa^+)^{\M}$. By induction on $n<\omega$ we construct a sequence $\langle (p_n, X_n):\,n<\omega\rangle$ such that $p_0=\emptyset_{\mathbb P}$ and for every $n<\omega$:
\begin{itemize}
\item $p_n\le_{n+1} p_{n+1}$, 
\item $p_{n+1} \forces``\name{h}(n) \in X_n"$ and
\item $|X_n|\le\kappa$.
\end{itemize}
We use (1) at successor stages. Let $q$ be the fusion of $\langle p_n:\,n<\omega\rangle$ and $X=\bigcup_{n<\omega} X_n$. By the definability of the forcing relation in $\M$ we have that $X\in \M$, and $X$ is a set of size $\le\kappa$. Therefore
$\alpha=\sup(X)<\kappa^+$, while $q\forces ``\sup(\ran(\name{h}))\le\alpha"$. This is in contradiction with the choice of $\name{h}$.
$\eop_{\ref{Miller-density}}$
\end{proof}

What remains is to show is that no bounded sequences are added to $\kappa$. In \S\ref{sec:with-filters}, we shall exhibit a whole family of forcings which can be used to solve the problem, including the forcing above. We note that for what was done so far, we could have chosen to demand the splittings at each stage $n$ to simply form an unbounded set in $\kappa_{n+1}$, but the argument about bounded sequences requires the splittings to come from a filter of sets.

\section{Generalised Perfect Set Forcing at a singular cardinal and the solution to (SP) for countable cofinality}\label{sec:with-filters} We use the same notation from \S\ref{sec:the-original-Prikry-modification}.
In this section we prove the case of Theorem \ref{the:main-result-unc} for singular cardinals of countable cofinality.

We start by modifying the forcing ${\mathbb P}_{\bar{\kappa}}$ to obtain a whole family of forcings that can be used. In addition to the sequence
$\bar{\kappa}$, we also choose a sequence $\bar{\FF}=\langle \FF_n:\,n<\omega\rangle$, where 
$\FF_0=\emptyset$ and 
\begin{itemize}
\item each $\FF_{n+1}$ is a $(<\kappa_n)$-complete filter on $\kappa_{n+1}$, 
\item $\FF_{n+1}$ does not contain any singletons and
\item $\FF_{n+1}$ does not contain any bounded subsets of $\kappa_n$. 
\end{itemize}

\begin{example}\label{ex:filters} We give some examples of filters $\FF_{n+1}$ which fit the above requirements.

{\noindent (1)} $\FF_{n+1}$ is simply the singleton $\{\kappa_{n+1}\}$. Taking all $\FF_{m+1}$ of this form is a special case of our argument which gives the forcing considered in \S\ref{sec:the-original-Prikry-modification}.

\smallskip
{\noindent (2)} $\FF_{n+1}$ is the filter of co-bounded subsets of $\kappa_{n+1}$.

\smallskip
{\noindent (3)} $\FF_{n+1}$ is a non-principal $(<\kappa_n)$-complete ultrafilter on $\kappa_{n+1}$. The latter choice carries a large cardinal strength: for uncountable cardinals  $\lambda<\theta$ the existence of a non-principal $(<\lambda)$-complete ultrafilter on $\theta$ implies that there is a measurable cardinal in the interval  $[\lambda, \theta]$. This is so since a modification of an argument of Ulam from \cite{Ulam} (see \cite[Lemma 10.2]{Jec03}) implies that the first cardinal that carries such an ultrafilter must be measurable. 
\end{example}

The forcing we define will be denoted
${\mathbb P}_{\bar{\kappa}} [\bar{\FF}]$. The conditions in ${\mathbb P}_{\bar{\kappa}} [\bar{\FF}]$ are defined similarly to those in
${\mathbb P}_{\bar{\kappa}}$, except that we modify what we mean by {\em fully-splitting}:

\begin{definition}\label{def:FSpl} (1) We say that a node $(n,t)$ in a perfect tree $T\subseteq T^{\ast}$ is  {\em fully-splitting for $\bar{\FF}$} if
\[
\{t'(n):\,(n,t)<_T (n+1,t')\}\in \FF_{n+1}.
\]

{\noindent (2)} The conditions, the order and the orders $\le_n$ are defined in ${\mathbb P}_{\bar{\kappa}} [\bar{\FF}]$ in the same way as in 
${\mathbb P}_{\bar{\kappa}}$, except that the requirement on fully-splitting is everywhere replaced by fully-splitting for $\bar{\FF}$.
$\eod_{\mbox{\tiny{Def.}}\ref{def:FSpl}}$
\end{definition}

\begin{claim}\label{obs:fusion+generic}
Lemma \ref{lem:generic}, Lemma \ref{Miller-density} and Lemma \ref{lem:fusion} hold with ${\mathbb P}_{\bar{\kappa}}$ replaced by ${\mathbb P}_{\bar{\kappa}} [\bar{\FF}]$. 
\end{claim}

\begin{proof} The proofs are very much the same as before, but let us comment on the density arguments of Lemma \ref{Miller-density}, using the same notation as in the proof of that lemma. For the density of the sets $\EE_g$ we need to know that $\bar{\FF}$ consists of non-principal filters, so $t' (o)$ can be chosen outside of the set $\{g(o)\}$. 

Regarding the density of the sets $\DD_{n,\zeta}$, we have assumed that $\FF_{o+1}$ does not contain bounded subsets of $<\kappa_m$.
Therefore the set of ordinals of the form $t'(o)$ must contain an ordinal $>\zeta$.
$\eop_{\ref{obs:fusion+generic}}$
\end{proof}

In the rest of the paper, we shall continue to use the notation $f=g[G]$ as in the above. From this point on we keep $\bar{\kappa}$ and
$\bar{\FF}$ fixed and use ${\mathbb P}$ to denote ${\mathbb P}_{\bar{\kappa}}[\bar{\FF}]$.

Next, in Lemma \ref{lem:no-bounded} we prove that 
${\mathbb P}_{\bar{\kappa}} [\bar{\FF}]$ does not add any bounded subsets to $\kappa$, so completing the countable cofinality case of the proof of Theorem  \ref{the:main-result-unc}.

We need to introduce certain auxiliary orders, which are not the same as the $\le_n$ we used in the fusion argument. 

\begin{definition}\label{def:more-orders} (1) 
For $p,q \in {\mathbb P}$ we write $\equiv_n$ for the equivalence relation which holds if $p, q$ have the same first $n$-levels. 


{\noindent (2)} For $p,q \in {\mathbb P}$, we write $p\preceq _n q$ if $p\le q$ and and $p \equiv_n q$.


$\eod_{\mbox{\tiny{Def.}} \ref{def:more-orders}}$
\end{definition}

\fbox{Wrong} The following Lemma is incorrect, which invalidates the main argument.

\begin{lemma}\label{lem:closure-of-len} For every $n<\omega$, the order $\preceq_n$ is $(<\kappa_{n})$-closed in ${\mathbb P}$. 
\end{lemma}

\begin{proof} Suppose that $\zeta^\ast<\kappa_{n}$ and $\bar{p}=\langle p_\zeta:\,\zeta<\zeta^\ast\rangle$ is increasing in $\preceq_n$. Then 
$q=\bigcap_{\zeta<\zeta^\ast} p_\zeta$ is a common extension of $\bar{p}$, because all $p_\zeta$ have the same first $n$ levels, and the splitting on higher levels come from filters which are at least $(<\kappa_{n})$-complete.
$\eop_{\ref{lem:closure-of-len}}$
\end{proof}

\begin{lemma}[{\bf Main Lemma for Theorem \ref{the:main-result-unc}, countable cofinality}]\label{lem:no-bounded} The forcing ${\mathbb P}_{\bar{\kappa}} [\bar{\FF}]$ does not add bounded subsets to $\kappa$.
\end{lemma}

\begin{proof} Suppose that $\name{X}$ is a $\mathbb P$-name for a subset of 
$\kappa_n$ for some $n$. Let $p\in P$ and $\alpha<\kappa_n$. We first construct $q=q_\alpha(p)$ such that 
\begin{itemize}
\item $p\preceq_{n+1} q$ and 
\item whenever $r, r'\ge q$ are such that $\lg(s[r]), \lg(s[r'])\ge n+1$ and $r,r'$ decide if $\alpha\in \name{X}$, then
\[
s[r]\rest (n+1)=s[r']\rest (n+1)\implies [r \forces ``\alpha\in \name{X}" \iff r' \forces ``\alpha\in \name{X}"].
\]
\end{itemize}

To construct $q$, first fix $(n+1,s)\in p$. The condition $p_{((n+1),s)}$ is an extension of $p$, so it forces that $\name{X}$ is a subset of $\kappa_n$.
Hence, it has an extension which decides if $\alpha\in \name{X}$. Such an extension is necessarily of the form $p_{(k,t)}$ for some $(k,t)\ge_p 
((n+1),s)$, so we choose some such extension and we call it $q^{(n+1,s)}$.
We choose such $q^{(n+1,s')}$ for every $s'$ such that $(n+1,s')\in p$. and define 
\[
q=\bigcup_{\{s':\, (n+1, s')\in p\}} q^{(n+1,s')}.
\]
We claim that $q$ is as required. Clearly $p \preceq_{n+1} q $ since $q$ is a subtree of $p$ and every $(n+1,s)\in p$ is contained in $q^{(n+1,s)}$, and $q^{(n+1,s)}\subseteq p$. 

Now suppose that $r,r'\ge q$ are conditions that decide if $\alpha\in \name{X}$ whose both stems are of length at least $n+1$, and such that $s[r]\rest (n+1)=s[r']
\rest (n+1)=(n+1,s)$. In particular, both $r$ and $r'$ are compatible with $q^{(n+1,s)}$ and therefore $r$ and $r'$ force that $\alpha\in {X}$ iff
$q^{(n+1,s)}$ does.

Having shown how to construct $q_\alpha(p)$ for a fixed $p$ and $\alpha$, we proceed to define by induction on $\alpha<\kappa_n$ a sequence 
$\langle p_\alpha:\,\alpha<\kappa_n\rangle$, as follows:
\begin{itemize}
\item $p_0=p$,
\item $p_{\alpha+1}=q_\alpha(p_\alpha)$,
\item $p_\delta=\bigcap_{\alpha<\delta} p_\alpha$ if $0<\delta<\kappa_n$ is a limit ordinal.
\end{itemize}
Notice that there is no problem in doing the limit steps, because the sequence is chosen to be $\prec_{n+1}$-increasing and by Lemma 
\ref{lem:closure-of-len}   this
order is $(<\kappa_{n+1})$-closed. For the same reason, defining $q^\ast=\bigcap_{\alpha<\kappa_n} p_\alpha$, we obtain a condition in $\mathbb P$.
This condition is an extension of $p$ and we claim that it forces that $\name{X}$ is an element of $\M$. Once we show this, the proof of the lemma is complete, since we have effectively shown that the set of conditions that force that $\name{X}$ is an element of $\M$, is dense.

So suppose that $r\ge q^\ast$ forces that $\name{X}$ is not an element of $\M$. Without loss of generality we can extend $r$ so that its stem has length
at least $n+1$. Let us denote its restriction to $n+1$ by $(n+1,s)$. Our construction was made so that 
\[
r \forces``\name{X}=Y=\{\alpha<\kappa_n:\,p_\alpha^{(n+1, s)}\forces ``\alpha\in \name{X}"\}
\]
and $Y$ is a set in $\M$. A contradiction.
$\eop_{\ref{lem:no-bounded}}$
\end{proof}

This finishes the proof of the countable cofinality case of Theorem \ref{the:main-result-unc}.

\section{Generalised Perfect Set Forcing at a regular cardinal-recall}\label{sec:perfect-regular} 
Let us recall the definition of the forcing $\mathbb P(\FF)$ introduced in \cite{Property-B}. The forcing 
${\mathbb P}_{\bar{\kappa}} [\bar{\FF}]$, which is our main interest in this paper, is an adaptation to a singular cardinal of the forcing $\mathbb P(\FF)$.
On the other hand, for uncountable $\kappa$, the forcing $\mathbb P(\FF)$ is the appropriate generalisation of the Perfect Set Forcing. At $\kappa=\aleph_0$, the direct predecessor of these definitions is Grigorieff forcing from \cite{zbMATH03513777}. We recall the relevant definitions in order to blend them with the method of the proof of the countable cofinality case, to pass to the arbitrary cofinality in Theorem \ref{the:main-result-unc}.

\begin{definition}\label{gSacks:def} 
Let $\kappa$ be an infinite regular cardinal and let $\FF$ be a $(<\kappa)$-complete filter on $\kappa$. We say that a subtree $p\subseteq {}^{{<\kappa}}\kappa$ is  {\em $\FF$-perfect} if 
\begin{itemize}
\item $\sup\{\dom(s):\,s\in p\}=\kappa$, and
\item for any $s\in p$, there is
$t\supseteq s$ such that $\{\alpha<\kappa:\, t\frown\alpha\in p\} \in \FF$.
\end{itemize}
We say that such $t$ as in the last property is $\FF$-{\em splitting}. 
For $s\in p$ we define its {\em splitting history} as
\[
{\rm deg}_p (s)=\{i < \dom(s):\,(\exists t\in p)\,[t\rest i =s\rest i \mbox{ and }t(i)\neq s(i)]\}.
\]

(1) The $\FF$-perfect tree forcing $\mathbb P=\mathbb P(\FF)$ is defined as follows:

The conditions in $\mathbb P$ are $\FF$-perfect subtrees $p$ of ${}^{{<\kappa}}\kappa$, satisfying the following additional conditions:
\begin{description}
\item [{\rm (o)}] For every $s\in p$, either $| \{\alpha<\kappa:\,s\frown \alpha\in p\} | = 1$, or $\{\alpha<\kappa:\,s\frown \alpha\in p\}\in \FF$.
\item [{\rm (a) [Closure of $p$]}] If $\delta<\kappa$ is a limit ordinal and $s\in {}^{\delta} \kappa$ is such that for all $\alpha<\delta$ we have $s\rest\alpha\in p$, then $s\in p$.
\item [{\rm (b) [Closure of splitting for every node]}] For every $s\in p$, the set ${\rm deg}_p (s)$ is a closed subset of $\dom(s)$.
\end{description}
The order is defined by $p\le q$ if $p\supseteq q$.

The {\em stem} of $p$,
denoted by $s[p]$ is defined by
\[
s[p]=\bigcup \{s\in p:\,\mbox{no }t\subset s \mbox{ is splitting and }s \mbox{ is splitting}\}.
\]
(Note that this is well-defined by (o) and that in fact $s[p]$ is a uniquely defined element of $p$ since there cannot be two incomparable $s$ satisfying the conditions in the definition of $s[p]$).

\smallskip

{\noindent (2)} For $\zeta<\kappa$ the order $\le_\zeta$ is defined by letting 
$p\le_\zeta q$ if 
\begin{itemize}
\item $p\le q$,
\item for all $s\in p$ with ${\rm otp}({\rm deg}_p(s))< \zeta$ we have $s\in q$.
\end{itemize}
\end{definition} 

The following were proved in \cite{Property-B} and \cite{Kanovei-style}:

\begin{fact}\label{fact:known}
(1) Every sequence of the form $\bar{p}=\langle p_\zeta:\,\zeta<\kappa\rangle$ such that $p_\zeta\le_\zeta p_{\zeta+1}$ for each $\zeta$ and $p_\delta=
\bigcap_{\zeta<\delta} p_\zeta$ for $\delta\in (0,\kappa)$ a limit ordinal, has a least upper bound which is $q=\bigcap_{\zeta<\kappa} p_\zeta$.
It in addition satisfies $p_\zeta\le_\zeta q$, for every $\zeta <\kappa$.

\smallskip

{\noindent (2)} Suppose that $\kappa>\aleph_0$.
The generic function $g=g_G:\,\kappa\to\kappa$ added by $\mathbb P(\FF)$ satisfies that for every $f:\,\kappa\to\kappa\in M$ there is a club $C=C_f$ of $\kappa$ in $M$
such that for every $\alpha\in C$ we have $f(\alpha)<g(\alpha)$.
\end{fact} 

A sequence $\bar{p}$ as in Fact \ref{fact:known}(2) is called a {\em fusion sequence}. 

\section{The proof of Theorem \ref{the:main-result-unc}, completing the solution to (SP)}\label{sec:unc-cof} In this section we are going to adapt the proof of 
the countable cofinality case to singular cardinals $\kappa$ of uncountable cofinality, and so obtain the complete proof of Theorem \ref{the:main-result-unc}

We start with a singular cardinal $\kappa$ with $\theta=\cf(\kappa)<\kappa$ arbitrary. Since we have already dealt with the countable cofinality, for notational purposes we assume $\theta>\aleph_0$. We also assume that $\kappa$ is a strong limit, which is necessary for the Mastering $\kappa$ argument in this case - see the explanation after the definition of $T^\ast$ below. Let us fix an increasing sequence $\bar{\kappa}=\langle \kappa_i:\,i<\theta\rangle$ of regular cardinals satisfying:
\begin{itemize}
\item $\kappa_0=\theta$,
\item $\kappa_{i+1}\ge 2^{\kappa_i}$ (this is possible, since $\kappa$ is a strong limit),
\item $\kappa_\delta\ge\sup_{i<\delta}\kappa_i$ for $\delta\in(0,\theta)$ a limit ordinal and
\item $\sup_{i<\theta} \kappa_i=\kappa$.
\end{itemize}
We also choose a sequence $\bar{\FF}=\langle \FF_{i+1}:\,i<\theta\rangle$, where each $\FF_{i+1}$ is a 
$(<\kappa_i)$-complete non-principal filter on $\kappa_{i+1}$ (might be just the singleton $\{\kappa_{i+1}\}$). The forcing we define will be denoted
${\mathbb P}_{\bar{\kappa}} [\bar{\FF}]$. 

We shall work with subtrees of the structure
\[
T^\ast=T^\ast_{\bar{\kappa}}=\bigcup_{i<\theta} (\{i\}\times {}^{i}\kappa_i),
\]
where the elements of $T^\ast$ are partially ordered by the following relation:
\[
(i, s)<^\ast (j,t)
\]
if $i<j$ and $s\subseteq t$ or $n=0$ and $s=\langle \rangle$ is the empty sequence. We introduce the concept of generalised perfect subtrees of 
$T^\ast$, which is similar to what we did in Definitions \ref{def:perfect} and \ref{def:FSpl}. However, we embed into the definition the requirements from Definition \ref{gSacks:def}(1).

In the argument for Mastering $\kappa$, similar to that of Lemma \ref{Miller-density} , and hence for the preservation of $\kappa^+$ we need that $|T^\ast|=\kappa$.
This was automatic in the case $\cf(\kappa)=\aleph_0$, but in general we need to assume that $\kappa$ is a strong limit.

\begin{definition}\label{def:perfect-unc} Suppose that $p\subseteq T^\ast$ is a tree.

(1) We say that a node $(i,s)\in p$ is {\em terminal} if $|\{t:\;(i,s)<_p(i, t)\}|= 0$, otherwise the node is {\em non-terminal}. The non-terminal node
$(i,s)$ is  {\em splitting in} $p$ if $|\{t:\;(i,s)<_p(i+1, t)\}|\ge 1$, otherwise the non-terminal node $(i,s)$ is a non-splitting node.
A splitting node $(i,s)$ is {\em fully-splitting for $\bar{\FF}$} if 
\[
\{t(i):\;(i,s)<_p(i+1, t)\}\in \FF_{i+1}.
\]
For $(i,s)\in p$ we define its {\em splitting history} as
\[
{\rm deg}_p ((i,s))=\{j < i:\,(j,s\rest j) \mbox{ is fully-splitting for }\bar{\FF}\}.
\]
{\noindent (2)} We say that a tree $p \subseteq T^\ast$ is {\em $(\bar{\kappa}, \bar{\FF})$-perfect} if
\begin{itemize}
\item $\sup\{i:\,(\exists s)(i,s)\in p\}=\theta$, 
\item for any $(i,s)\in p$, there is
$(j,t)>_p (i, s)$ which is fully-splitting for $\bar{\FF}$,
\item (a) [Closure of $p$] If $\delta<\theta$ is a limit ordinal and $s\in {}^{\delta} \kappa$ is such that for all $i<\delta$ we have $(i,s\rest i)\in p$, then 
$(\delta,s)\in p$, and
\item (b) [Closure of splitting for every node] For every $(i,s)\in p$, the set ${\rm deg}_p ((i,s))$ is a closed subset of $i$.
\end{itemize}

{\noindent  (3)} The conditions in $\mathbb P=\mathbb P_{\bar{\kappa}}[\FF]$ are $(\bar{\kappa}, \bar{\FF})$-perfect subtrees $p$ of $T^\ast$, which we order by reverse inclusion.

{\noindent (4)} For $p,p' \in \mathbb P$ and $\zeta<\theta$, we let $p\le_\zeta p'$ if $p\le p'$ and for every $(i,s)\in p$,
\[
\deg_p((i,s))\le \zeta\implies (i,s)\in p'.
\]
$\eod_{\mbox{\tiny{Def.}} \ref{def:perfect-unc}}$
\end{definition}

The proof of Theorem \ref{the:main-result-unc} now proceeds exactly as in the case of countable cofinality, but using the version of the forcing $\mathbb P_{\bar{\kappa}}[\FF]$
as presented here.

\section*{Research Support Acknowledgements}  
The author used the AI product GROK by X-AI as a research assistant for rough historical facts, finding references and \LaTeX questions. All mathematical developments are due to the author unassisted by the AI. The author thanks the British Library in London, the Bibliothèque Nationale de France, l'Institut Henri Poincare in Paris, la Bibliothèque Royale de Brussels, Math-Net.Ru, Association for Symbolic Logic and the London Mathematical Society for the memberships, providing her access to the scientific resources necessary for her research. She especially thanks the Steklov Institute of Mathematics in Moscow for their support during her visit for the week of 16th November 2026.

The mathematical support of being able to discuss with Vladimir Kanovei on this and related subjects is priceless.

\bibliography{../bibliomaster}
\bibliographystyle{plainurl}

\end{document}